\documentclass{article}
\usepackage{arxiv}

\usepackage[utf8]{inputenc} 
\usepackage[T1]{fontenc}    
\usepackage{hyperref}       
\usepackage{url}            
\usepackage{booktabs}       
\usepackage{amsfonts}       
\usepackage{nicefrac}       
\usepackage{microtype}      
\usepackage{lipsum}		
\usepackage{graphicx}
\usepackage{doi}
\usepackage{diagbox}
\usepackage{float}
\usepackage{multicol}
\usepackage{caption}
\usepackage{subfig}
\usepackage[symbol]{footmisc}

\usepackage[square,sort,comma,numbers]{natbib}

\usepackage[dvipsnames]{xcolor}
\definecolor{blue}{RGB}{0,56,118}
\definecolor{red}{RGB}{166,25,46}
\definecolor{gold}{RGB}{109,93,5}

\title{\Large How a Losing Team like the Canadiens can Steal a Stanley Cup: 
A Quantitative Intransitive Hockey Analysis}

\author{
C. J. Barrett\textsuperscript{*}, S. Koumarianos, O. Mermut,
York University (Toronto); McGill University (Montreal), Canada. 
}

\date{\textbf{submitted:} 17 June 2021, to:  arXiv.org}

\renewcommand{\headeright}{}
\renewcommand{\undertitle}{}
\renewcommand{\shorttitle}{}

\begin{document}

\maketitle

\textit{Dedicated to the Canadian academics who volunteered their mathematical expertise to epidemiological modeling 2020–2021. On the occasion of Canada's latest hockey playoff effort towards bringing Lord Stanley’s Cup back home for the first time since 1993, by pulling off one of the most unlikely underdog upsets in professional hockey history.\\}

\begin{abstract}
We present here a simple mathematical model that provides a successful strategy, quantitatively, to ending the continued championship futility experienced by Canadian Hockey Teams.  Competitive Intransitivity is used here as a simple predictive framework to capture how investing strategically, under a uniform salary cap, in just 3 independently variable aspects of the sport (such as \textit{Offence}, \textit{\textit{Defence}}, and a \textit{Goaltender}), by just 3 Hockey Teams applying differing salary priorities (such as \textcolor{red}{Montreal}, \textcolor{gold}{Boston}, and \textcolor{blue}{New York}), can lead to rich and perhaps surprisingly unexpected outcomes in play, similar to rolling intransitive dice together in a series of head-to-head games. A possibly fortunate conclusion of this analysis is the prediction that for any Team's chosen strategy (such as \textcolor{blue}{New York}'s), a counter strategy within the same salary cap can be adopted by a playoff opponent (such as \textcolor{red}{Montreal}) which will prove victorious over a long playoff series, enabling a pathway to end prolonged championship futility.
\end{abstract}

\section{Assumptions of this Model}
We construct here a simple description of Hockey Playoffs as between just 3 Teams (such as: \textcolor{red}{Montreal}, \textcolor{gold}{Boston}, and \textcolor{blue}{New York}), where each Team possesses different strengths in just 3 independent competitive variables (such as \textit{Offence}, \textit{\textit{Defence}}, and a \textit{Goalie}), represented by different whole numbers (such as \$ millions), summing to the same total (a 'salary cap' such as \$6 million  /Team). Such '\textit{Goalie}-centred', 'balanced' and '\textit{Offence}-\textit{\textit{Defence}}' spending could be represented, for example as:
{
\begin{multicols}{2}
\begin{center}
    \textbf{}\\\textbf{}\\
    \renewcommand{\arraystretch}{1.3}
   \begin{tabular}{cccc}
    &\textit{\textit{Offence}} (\$M) &\textit{\textit{\textit{Defence}}} (\$M) &\textit{\textit{Goalie}} (\$M) \\
    \textcolor{red}{\textcolor{red}{\textcolor{red}{Montreal}}} & 1 & 1  & 4  \\
    \textcolor{gold}{\textcolor{gold}{Boston}} & 2 & 2  & 2  \\
    \textcolor{blue}{\textcolor{blue}{New York}} & 3 & 3  & 0  
\end{tabular} 
\end{center}
\columnbreak
 \begin{center}
 \textbf{}\\
      \includegraphics[scale=0.3]{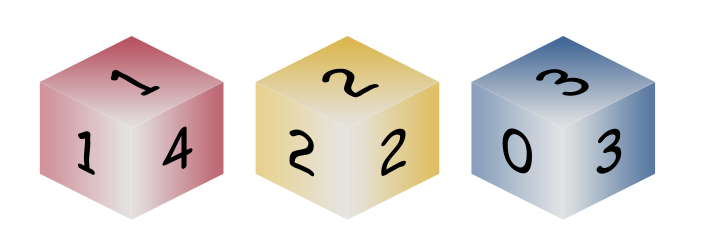}  
 \end{center}
\end{multicols}
}

The Model is run by assuming that each pair of Teams plays each other over a long series (approaching $\infty$), and that the winner of that series is the Team who wins the most 'head-to-head matchups' of these 9 possible combinations of competitive variables, similar to rolling differing dice against each other many times, to see which die 'wins'.\\

Which strategy is best? \textit{ i.e.: }is it better for \textcolor{blue}{New York} to spend so much on \textit{Offence} and \textit{\textit{Defence}}, or for \textcolor{red}{Montreal} to concentrate resources in their \textit{Goalie}, or does \textcolor{gold}{Boston} end up victorious with balanced spending? 

\newpage
\section{Results from the Model}
\label{sec:headings}

 The 9 independent ‘head-to-head matchups’ between each pair of Teams facing each other in a Playoff Series, might be easily visualized as rolling 3 different coloured dice, representing the 3 Teams’ weighting strategies in \textit{Offence}, \textit{\textit{Defence}}, and \textit{Goalie} variables (repeating the same 3 numbers on the backside of each 6-sided die).

\textbf{Playoff Series Winners} can be presented by charting results of the 9 possible match-ups, then declaring as winner the Team who out-rolls their opponent in the majority of the 9 possible combinations, for \textit{e,g,:} 

\begin{multicols}{2}
    \begin{center}
    \renewcommand{\arraystretch}{1.3}
    \begin{tabular}{|c|ccc|}
        \hline
        \diagbox{\textcolor{red}{MTL}}{\textcolor{gold}{BOS}} & 2 & 2 & 2\\
        \hline        
        1 & \textcolor{gold}{BOS} & \textcolor{gold}{BOS} & \textcolor{gold}{BOS} \\
        1 & \textcolor{gold}{BOS} & \textcolor{gold}{BOS} & \textcolor{gold}{BOS} \\ 
        4 & \textcolor{red}{MTL} & \textcolor{red}{MTL} & \textcolor{red}{MTL} \\
        \hline
    \end{tabular}
    \end{center}
        \begin{center}
        \renewcommand{\arraystretch}{1.3}
        \begin{tabular}{|c|ccc|}
            \hline
            \diagbox{\textcolor{gold}{BOS}}{\textcolor{blue}{NY}} & 3 & 3 & 0\\
            \hline
            2 & \textcolor{blue}{NY} & \textcolor{blue}{NY} & \textcolor{gold}{BOS}\\
            2 & \textcolor{blue}{NY} & \textcolor{blue}{NY} & \textcolor{gold}{BOS}\\
            2 & \textcolor{blue}{NY} & \textcolor{blue}{NY} & \textcolor{gold}{BOS}\\
            \hline
        \end{tabular}
        \end{center}
\end{multicols}

Where (left) in a matchup with \textit{Goalie}-heavy \textcolor{red}{Montreal}, a balanced \textcolor{gold}{Boston} Team would be expected to prevail eventually, ‘winning’ 6 of the possible 9 total matchups of Team strength.  Similarly (right), \textcolor{gold}{Boston} then playing an \textit{Offence}-\textit{\textit{Defence}} oriented \textcolor{blue}{New York} Team would be expected to be defeated, again in 6 out of 9 possible matchups. Since \textcolor{blue}{New York} triumphs over \textcolor{gold}{Boston}, after \textcolor{gold}{Boston} has clearly vanquished \textcolor{red}{Montreal}, one might be tempted to assume that a \textcolor{blue}{New York} \textit{vs} \textcolor{red}{Montreal} final would be as easily predictable as: \textcolor{blue}{NY} > \textcolor{gold}{BOS} > \textcolor{red}{MTL}, so therefore: \textcolor{blue}{NY} > \textcolor{red}{MTL} ?

\textbf{A Possibly Unexpected Final Outcome} can be confirmed by the match-up chart between \textcolor{red}{Montreal} and \textcolor{blue}{New York} (below), where examining results of the 9 combinations reveals that in only 4 of 9 matchups does \textcolor{blue}{New York} prevail, yet underdog \textcolor{red}{Montreal} emerges victorious overall, winning 5 of 9 matchups, and thus defeating the \textcolor{blue}{New York} Team.  Such possibly surprising yet fortuitous final outcomes can be described as ‘intransitive’, with much written elsewhere about such potentially unexpected and apparently undeserving outcomes, using many variations of such intransitive dice. 

\begin{center}
    \renewcommand{\arraystretch}{1.3}
    \begin{tabular}{|c|ccc|}
        \hline
        \diagbox{\textcolor{red}{MTL}}{\textcolor{blue}{NY}}&3&3&0\\
        \hline
        1 & \textcolor{blue}{NY} & \textcolor{blue}{NY} & \textcolor{red}{MTL}\\
        1 & \textcolor{blue}{NY} & \textcolor{blue}{NY} & \textcolor{red}{MTL}\\
        4 & \textcolor{red}{MTL} & \textcolor{red}{MTL} & \textcolor{red}{MTL} \\
        \hline
    \end{tabular}
\end{center}

\section{Conclusions}
 It is demonstrated here by this Model that for any distribution of funding adopted by any Team (for example: \textcolor{gold}{Boston} or \textcolor{blue}{New York}) under a uniform salary cap, a superior winning distribution of the same resources can be adopted by their opponent (such as: \textcolor{red}{Montreal}) to ensure final victory, and a welcomed end to prolonged championship futility. 

\section{Acknowledgements \& Competing Interests}
Profs. C.B. and O.M. are grateful to NSERC Canada for research support. The Canadian Broadcasting Corporation is thanked for open broadcast of \textit{‘Hockey Night in Canada’} during which this collaborative paper was written, and Moosehead Breweries Limited (St. John, NB). The authors declare no competing interests, aside from traditional geographical hockey allegiances.

\small{
\bibliographystyle{unsrt}

\begin{thebibliography}{1}

\bibitem{wiki}
{Wikipedia contributors}.
\newblock Intransitive dice --- {W}ikipedia{,} the free encyclopedia, 2021.
\newblock [Online; accessed 24-May-2021].

\bibitem{gardner1970paradox}
Martin Gardner.
\newblock The paradox of nontransitive dice and the elusive principle of indifference.
\newblock {\em Scientific American}, 223(6):110, 1970.

\bibitem{leonard2010neumann}
Robert Leonard.
\newblock {\em Von Neumann, Morgenstern, and the creation of game theory: From
  chess to social science, 1900--1960}.
\newblock Cambridge University Press, 2010.

\bibitem{ekhad2017treatise}
Shalosh~B Ekhad and Doron Zeilberger.
\newblock A treatise on sucker's bets.
\newblock {\em arXiv preprint arXiv:1710.10344}, 2017.

\end{thebibliography}

}
\end{document}